\documentclass{elsart}

\usepackage{amsmath,amssymb,amsfonts,amssymb,amsopn,amscd}
\usepackage{bbm,wasysym,stmaryrd}
\usepackage[english,francais]{babel}
\usepackage{color}
\usepackage{hyperref}

\newtheorem{theorem}{Theorem}[section]
\newtheorem{lemma}[theorem]{Lemma}
\newtheorem{e-proposition}[theorem]{Proposition}

\newtheorem{e-definition}[theorem]{Definition\rm}

\newtheorem{theoreme}{Th\'eor\`eme}[section]

\newtheorem{proposition}[theoreme]{Proposition}

\newtheorem{definition}[theoreme]{D\'efinition\rm}

\newcommand{\olivier}[1]{{#1}}
\newcommand{\R}{\mathbbm{R}}
\newcommand{\mN}{\mathcal{N}}
\newcommand{\Exp}{\mathbb{E}} 
\newcommand{\Expt}[1]{\mathbb{E}\left [ #1 \right]} 
\newcommand{\J}{J}
\newcommand{\Z}{\mathbbm{Z}}
\newcommand{\T}{\mathcal{T}}  
\newcommand{\mS}{\mathcal{S}}  
\newcommand{\M}{\mathcal{M}}

\newcommand{\mmeas}{\mathcal{M}_{\mS}}
\newcommand{\B}{\mathcal{B}}
\newcommand{\mM}{\mathcal{M}}
\newcommand{\mQ}{\mathcal{Q}}
\newcommand{\undt}[1]{\underline{#1}}
\newcommand{\mone}{1}

\def\og{\leavevmode\raise.3ex\hbox{$\scriptscriptstyle\langle\!\langle$~}}
\def\fg{\leavevmode\raise.3ex\hbox{~$\!\scriptscriptstyle\,\rangle\!\rangle$}}

\journal{the Acad\'emie des sciences}
\begin{document}
\centerline{}
\begin{frontmatter}


\selectlanguage{english}
\title{Asymptotic description of stochastic neural networks. I - existence of a Large Deviation Principle}


\selectlanguage{english}
\author[authorlabel1]{Olivier Faugeras},
\ead{firstname.name@inria.fr}
\author[authorlabel1]{James Maclaurin}

\address[authorlabel1]{Inria Sophia-Antipolis M\'editerran\'ee\\
NeuroMathComp Group}



\begin{abstract}
\selectlanguage{english}
We study the asymptotic law of a network of interacting neurons when the number of neurons becomes infinite. 
The dynamics of the neurons is described by a set of stochastic differential equations in discrete time.
The neurons interact through the
synaptic weights which are Gaussian {\emph correlated} random variables. We describe the 
asymptotic law of the network when the number of neurons goes to
infinity. Unlike previous works which made the biologically unrealistic assumption that the weights were i.i.d. random variables, we assume that they are correlated. We introduce the process-level empirical measure of the trajectories of the
solutions to the equations of the finite network of  neurons  and the averaged law (with respect to
the synaptic weights)  of the trajectories of the solutions to the equations of the network of neurons. The result (theorem \ref{theo:LDP} below) is that the image law through the empirical measure satisfies a large deviation principle  with a good rate function. We provide an analytical expression of this rate function in terms of the spectral representation of certain Gaussian processes.
\vskip 0.5\baselineskip

\selectlanguage{francais}
\noindent{\bf R\'esum\'e} \vskip 0.5\baselineskip \noindent
\begin{center}{\bf Description asymptotique de r\'eseaux de neurones stochastiques. I - existence d'un principe de grandes d\'eviation }\end{center}
\noindent
Nous consid\'erons un r\'eseau de neurones d\'ecrit par un syst\`eme d'\'equations diff\'erentielles stochastiques en temps discret. Les neurones interagissent au travers de poids synaptiques qui sont des variables al\'eatoires gaussiennes {\emph corr\'el\'ees}. Nous caract\'erisons la loi asymptotique de ce r\'eseau lorsque le nombre de neurones tend vers l'infini. Tous les travaux pr\'ec\'edents faisaient l'hypoth\`ese, irr\'ealiste du point de vue de la biologie, de poids ind\'ependants. Nous introduisons la mesure empirique sur l'espace des trajectoires solutions des \'equations du r\'eseau de neurones de taille finie et la loi moyenn\'ee (par rapport aux poids synaptiques) des trajectoires de ces solutions. Le r\'esultat (th\'eor\`eme \ref{theo:LDP} ci-dessous) est que l'image de cette loi par la mesure empirique satisfait un principe de grandes d\'eviations avec une bonne fonction de taux dont nous donnons une expression analytique en fonction de la repr\'esentation spectrale de certains processus gaussiens.
\end{abstract}
\end{frontmatter}

\selectlanguage{francais}
\section*{Version fran\c{c}aise abr\'eg\'ee}
Nous consid\'erons le probl\`eme de d\'ecrire la dynamique asymptotique d'un ensemble de $2n+1$ neurones lorsque ce nombre tend vers l'infini. Ce probl\`eme est motiv\'e par un d\'esir de parcimonie dans la description, par celui de rendre compte de l'apparition de ph\'enom\`enes \'emergents, ainsi que par celui de comprendre les effets de taille finie. Nous consid\'erons donc un r\'eseau de $2n+1$ neurones interconnect\'es dont la dynamique commune (en temps discret) ob\'eit aux \'equations stochastiques \eqref{eq:U}. Dans celles-ci apparaissent les poids synaptiques ou coefficients de couplage not\'es $J_{ij}^n$ qui sont des variables al\'eatoires gaussiennes corr\'el\'ees. Pour r\'epondre \`a la question pos\'ee nous consid\'erons la loi, not\'ee $Q^{V_n}$, de la solution \`a \eqref{eq:U} moyenn\'ee par rapport aux poids synaptiques ou plus pr\'ecis\'ement l'image $\Pi^n$ de cette loi par la mesure empirique \eqref{defn:hatmun}. Nous montrons dans le th\'eor\`eme \ref{theo:LDP} que cette loi satisfait un principe de grande d\'eviations avec une bonne fonction de taux $H$ dont nous donnons une expression analytique dans la d\'efinition \ref{defn:H} et les \'equations \eqref{eq:Gamma1infinity} et \eqref{eq:Gamma2}. Ce travail g\'en\'eralise au cas des poids synaptiques corr\'el\'es celui d'auteurs comme Sompolinsky \cite{sompolinsky-crisanti-etal:88} et Moynot et Samuelides \cite{moynot-samuelides:02} qui ont consid\'er\'e le cas de poids synaptiques ind\'ependants. Dans ce cas, plus simple d'un point de vue math\'ematique, mais beaucoup moins r\'ealiste d'un point de vue biologique, on observe le ph\'enom\`ene de propagation du chaos. Nous montrons dans un second article \cite{faugeras-maclaurin:14b} que la bonne fonction de taux a un minimum unique que nous caract\'erisons compl\`etement. La propagation du chaos n'a pas lieu mais la repr\'esentation est parcimonieuse dans un sens d\'efini dans \cite{faugeras-maclaurin:14b}.

\selectlanguage{english} \vspace{-1cm}
\section{Introduction}
\vspace{-0.5cm}
\subsection{Neural networks}
\vspace{-0.5cm}

Our goal is to study the asymptotic behaviour and large deviations of a network of
interacting neurons when the number of neurons becomes
infinite. A more detailed exposition of this work, with proofs, may be found in \cite{faugeras-maclaurin:14d}.

Sompolinsky 
succesfully explored this particular topic
\cite{sompolinsky-crisanti-etal:88} for fully connected networks of
neurons. In his
study of the continuous time dynamics of networks of rate neurons,
Sompolinsky and his colleagues assumed that the synaptic
weights, were i.i.d. random variables 
with zero mean Gaussian laws. The main result they obtained (using the local chaos hypothesis) under the
previous hypotheses is that the averaged law of the neurons
dynamics is chaotic in the sense that the averaged law of a finite
number of neurons converges to a product measure as the system gets very large.

The next efforts in the direction of understanding the
averaged law of neurons are those of Cessac, Moynot and Samuelides
\cite{cessac:95,moynot:99,moynot-samuelides:02,cessac-samuelides:07,samuelides-cessac:07}. From the
technical viewpoint, the study of the collective dynamics is done in discrete
time. Moynot and Samuelides obtained a large deviation principle and were able to describe in detail the
limit averaged law that had been obtained by Cessac using the local chaos hypothesis and to  prove rigorously the propagation of chaos
property. 

One of the next outstanding challenges is to incorporate in the network model the fact
that the synaptic weights are not independent and in effect, according to experimentalists, often
highly correlated. Our problem thus resembles that of a random walk in a mixing random environment \cite{sznitman-zerner:99,rassoul-agha:03}.

The problem whose solution we announce in this paper  and in \cite{faugeras-maclaurin:14b} is the following. Given a
completely connected network of neurons in which the
synaptic weights are Gaussian {\emph correlated} random variables, can we describe the 
asymptotic law of the network when the number of neurons goes to
infinity? \vspace{-1.0cm}

\subsection{Mathematical framework}
\vspace{-0.5cm}
For some positive integer $n> 0$, we let $V_n = \lbrace j\in\Z: |j| \leq n \rbrace$, and $|V_n|=2n+1$. The finite-size neural network below is indexed by points in $V_n$. We work in discrete time, over times $t\in \lbrace 0,1,\ldots,T\rbrace$, for some positive integer $T$. The state variable for each neuron is in $\R$, and the path space is $\T = \R^{T+1}$. We equip $\T$ with the Euclidean topology, $\T^{\Z}$ with the cylindrical topology, and denote the Borelian $\sigma$-algebra generated by this topology by $\B(\T^{\Z})$.

The equation describing the time variation of the membrane potential
$U^j$ of the $j$th neuron writes
\begin{equation}\label{eq:U}
U^j_{t}=\gamma U^j_{t-1}+\sum_{i\in V_n} J_{ji}^{n}
f(U^i_{t-1})+\theta^j+B^j_{t-1},  \quad U^j_0=u^j_0, \quad j \in V_n,\, t=1,\ldots,T
\end{equation}
\vspace{-0.5cm}

$f: \R \to ]0,\,1[$ is a monotonically increasing Lipschitz continuous bijection.
$\gamma$ is in $[0,1)$ and determines the time
scale of the intrinsic dynamics of the
neurons.
The $B^j_t$s are i.i.d. Gaussian random variables distributed as
$\mathcal{N}_1(0,\sigma^2)$\footnote{We note $\mN_p(m,\Sigma)$ the law of the $p$-dimensional
Gaussian variable with mean $m$ and covariance matrix $\Sigma$.}. They represent the fluctuations of the neurons' membrane potentials. The $\theta^j$s are i.i.d. as $\mathcal{N}_1(\bar{\theta},\theta^2)$. The are independent of the $B^i_t$s and represent the current injected in the neurons. The $u^j_0$s are i.i.d. random variables each governed by the law $\mu_I$.

The $J_{ij}^{n}$s are the synaptic weights. $J_{ij}^{n}$ represents the
strength with which the `presynaptic' neuron $j$ influences the
`postsynaptic' neuron $i$. They arise from a stationary Gaussian random field specified by its mean and covariance function  
\[
\Exp[J^n_{ij}]=\frac{\bar{J}}{|V_n|} \quad, cov(J_{ij}^nJ_{kl}^{n})=\frac{1}{|V_n|}\Lambda\left((k-i) \text{ mod } V_n,(l-j)\text{ mod }V_n\right),
\]
\vspace{-0.5cm}

$\Lambda$ is  positive definite, let $\tilde{\Lambda}$ be the corresponding (positive) Fourier transform. We make the technical assumption that the summation over both indices of the series $(\Lambda(i,j))_{i,j \in \Z}$ is absolutely convergent to $\Lambda^{sum}>0$.

We note $\J^n$ the $|V_n| \times |V_n|$ matrix of the synaptic weights, $\J^n=(J_{ij}^n)_{i,j \in V_n}.$

The process $(Y^j)$ defined by
\[
Y^j_t=\gamma Y^j_{t-1}+\bar{\theta}+B^j_{t-1},\quad j \in V_n,\quad t=1,\cdots T, \quad Y^j_0=u^j_0
\]
is stationary and independent. The law of each $Y^j$ is easily found to be given by
\[
P=(\mN_{T}(0_{T},\sigma^2{\rm Id}_{T}) \otimes \mu_I) \circ \Psi ,
\]
where $\Psi:\T \rightarrow \T$ is the following affine bijection. The joint law of $(Y^k)$ (for $k\in V_n$) is written as $P^{\otimes V_n}$, and the joint law of all $(Y^j)$ is written as $P^{\Z}$. Writing $v = \Psi(u)$, we define
\begin{equation}\label{eq:Psi}
\left\{
\begin{array}{lcl}
v_0 &=& \Psi_0(u)=u_0\\
v_s  &=&  \Psi_s(u)=u_s - \gamma u_{s-1}-\bar{\theta}\quad s=1,\cdots, T.
\end{array}
\right.
\end{equation}
We extend $\Psi$ to a mapping $\T^\Z \to \T^\Z$ componentwise. We now introduce some more notation.

For some topological space $\Omega$ equipped with its Borelian $\sigma$-algebra $\B(\Omega)$, we denote the set of all probability measures by $\M(\Omega)$. We equip $\M(\Omega)$ with the topology of weak convergence. For some $\mu\in\M(\T^{\Z})$ governing a process $(X^j)_{j\in\Z}$,  we let $\mu^{V_n}\in\M(\T^{V_n})$ denote the marginal governing $(X^j)_{j\in V_n}$. For some $\mu\in\M(\T^{\Z})$ governing a process $(X^j)_{j\in\Z}$,  we let $\mu^{V_n}\in \M(\T^{V_n})$ denote the marginal governing $(X^j)_{j\in V_n}$. For some $X \in \T$ and $0 \leq a \leq b \leq T$, $X_{a,b}$ denotes the $b-a+1$-dimensional subvector of $X$. We let $\mu_{a,b} \in \M(\T^\Z_{a,b})$ denote the marginal governing $(X^j_{a,b})_{j \in \Z}$. For some $j\in \Z$, let the shift operator $\mS^j:\T^{\Z}\to \T^{\Z}$ be $S(\omega)^k = \omega^{j+k}$. We let $\M_\mS(\T^{\Z})$ be the set of all stationary probability measures $\mu$ on $(\T^{\Z},\B(\T^{\Z}))$ such that for all $j\in\Z$, $\mu\circ (\mS^j)^{-1} = \mu$.

\begin{definition}\label{def:mubarbar}
For each measure $\mu \in \M(\T^{V_n})$ or $\mmeas(\T^\Z)$ we define $\undt{\mu}$ to be $\mu \circ \Psi^{-1}$. 
\end{definition}

We next introduce the following definitions.
\begin{definition}\label{def:E2}
Let  $\mathcal{E}_2$ be the subset  of $\mathcal{M}_\mS(\T^{\Z})$ defined by
\[
\mathcal{E}_2=\{ \mu \in \mathcal{M}_\mS(\T^{\Z})\, | \,
\Exp^{\undt{\mu}_{1,T}}[\|v^0\|^2] < \infty\}.
\]
\end{definition}
Let $p_n:\T^{V_n} \to \T^{\Z}$ be such that $p_n(\omega)^k = \omega^{k\mod V_n}$. Here, and throughout the paper, we take $k\mod V_n$ to be the element $l\in V_n$ such that  $l = k \mod |V_n|$.  Define the process-level empirical measure $\hat{\mu}_n: \T^{V_n} \to \M_\mS\left(\T^{\Z}\right)$ as
\begin{equation}
\hat{\mu}_n(\omega) = \frac{1}{|V_n|}\sum_{k\in V_n} \delta_{S^k p_n(\omega)}.\label{defn:hatmun}
\end{equation}

We define the process-level entropy to be, for $\mu\in\mathcal{M}_\mS(\T^{\Z})$
\begin{equation*}
I^{(3)}(\mu , P^\Z) = \lim_{n\to\infty} \frac{1}{|V_n|}I^{(2)}\left(\mu^{V_n},P^{\otimes V_n}\right).
\end{equation*}
If $\mu\notin \mathcal{E}_2$, then $I^{(3)}(\mu,P^\Z) = \infty$. Here $I^{(2)}$ is the \textit{relative entropy}. For further discussion, a definition of $I^{(2)}$ and a proof that $I^{(3)}$ is well-defined, see \cite{deuschel-stroock-etal:91}.

We note $Q^{V_n}(\J^n)$ the element of $\M(\T^{V_n})$ which is the law of the solution to \eqref{eq:U} conditioned on $\J^n$. We let $Q^{V_n} = \mathbb{E}^{\J}[Q^{V_n}(\J^n)]$ be the law averaged with respect to the weights. The reason for this is that we want to study the empirical measure $\hat{\mu}_n$ on path space. There is no reason for this to be a simple problem since for a fixed interaction $J^n$, the variables $(U^j)_{j \in V_n}$ are not exchangeable. So we first study the law of $\hat{\mu}_n$ averaged over the interactions. 

Finally we introduce the image laws in terms of which the principal results of this paper are formulated.
\begin{definition}\label{def:PiNQN}
Let $\Pi^n$ and $R^n$ in $\mathcal{M}(\mathcal{M}_\mS(\T^\Z))$ be the image laws of $Q^{V_n}$ and $P^{\otimes V_n}$ through the function
$\hat{\mu}_n: \T^{V_n} \to \mM_\mS(\T^\Z)$ defined by \eqref{defn:hatmun}:
\[
\Pi^n=Q^{V_n} \circ \hat{\mu}_n^{-1} \quad R^n=P^{\otimes V_n} \circ \hat{\mu}_n^{-1}
\]
\end{definition}
\vspace{-0.75cm}
\section{The good rate function}
\vspace{-0.75cm}
We obtain an LDP for the process with
correlations ($\Pi^n$) via the (simpler) process without correlations
($R^{n}$). To do this we obtain an expression for the Radon-Nikodym derivative of $\Pi^n$ with respect to $R^{n}$.  This is done in propositions \ref{prop:RNderiv1} and \ref{prop:radon-nikodym}. In equation \eqref{eq:dQNdPN} there appear certain Gaussian random variables defined from the right handside of the equations of the neuronal dynamics \eqref{eq:U}. Applying the Gaussian calculus to this expression we obtain equation \eqref{eq:RN1} which expresses the Radon-Nikodym derivative as a function (depending on $n$) of the empirical measure \eqref{defn:hatmun}. Using the fact that this function is measurable we obtain equation \eqref{eq:RN2}. This equation is essential in a) finding the expression for the function $\Gamma$  that appears in the rate function $H$ of definition \ref{defn:H}, 
b) proving the lower-bound for $\Pi^n$ on the open sets, c) proving that the sequence $(\Pi^n)$ is exponentially tight, and d) proving  the upper-bound on the compact sets.

The key idea is to associate to every stationary measure $\mu$ a certain stationary Gaussian process $G^\mu$, or equivalently a certain Gaussian measure defined by its mean $c^\mu$ and its covariance operator $K^\mu$. This allows us to write the Radon-Nikodym derivative as a function of the empirical measure, through writing is as a function of $G^{\hat{\mu}_n}$.

Given $\mu$ in $\mathcal{M}_\mS(\T^\Z)$ we define a stationary Gaussian process $G^\mu$, governed by a measure $\mQ^\mu \in \mathcal{M}_\mS(\T_{1,T}^\Z)$. For all $i$ the mean of $G^{\mu,i}_t$ is given by $c^\mu_t$, where
\begin{equation}
c^{\mu}_t = \bar{J} \int_{\T^\Z} f(u^i_{t-1})
d\mu(u),\, t=1,\cdots,T\, , i \in \Z,\label{eq:cmumean}
\end{equation}\vspace{-0.75cm}

The covariance between the Gaussian vectors $G^{\mu,i}$ and $G^{\mu,i+k}$ is defined to be\footnote{We note ${}^\dagger$ the transpose of a vector or matrix.}
\begin{equation}\label{eq:Kimuinfinite}
K^{\mu,k}=\theta^2 \delta_k \mone_{T} {}^\dagger  \mone_{T}+
\sum_{l=-\infty}^{\infty} \Lambda(k,l) M^{\mu,l},
\end{equation}\vspace{-0.75cm}

where $1_T$ is the $T$-dimensional vector whose coordinates are all equal to 1 and
\begin{equation}\label{eq:Mmuinfinite}
M^{\mu,k}_{st} = \int_{\T^\Z} f(u^0_{s-1}) f(u^k_{t-1}) d\mu(u),
\end{equation}\vspace{-0.75cm}

The above integrals are well-defined because of the definition of $f$ and the fact that the series in (\ref{eq:Kimuinfinite}) is convergent (since the series
$(\Lambda(k,l))_{k,\,l \in \Z}$ is absolutely convergent 
and the elements of
$M^{\mu,l}$ are bounded by $1$ for all $l \in \Z$).  We note $\mQ^\mu_{[n]}$ the law of the $|V_n|$-dimensional Gaussian defined by restricting the sum in \eqref{eq:Kimuinfinite} to $l \in V_n$.

These definitions imply the existence of a Hermitian-valued spectral representation for the sequence $M^{\mu,k}$ (resp. $K^{\mu,k}$) noted $\tilde{M}^\mu$ (resp. $\tilde{K}^\mu$) which satisfies
\[
\tilde{K}^\mu(\theta)=\theta^2\mone_T\,{ }^\dagger \mone_T+
\frac{1}{2\pi}\int_{-\pi}^{\pi} \tilde{\Lambda}(\theta,-\varphi) \tilde{M}^\mu(d\varphi).
\]
\vspace{-0.75cm}

This allows us to define the spectral representation
\begin{equation}\label{eq:Atildetheta}
\tilde{A}^{\mu}(\theta)=\tilde{K}^{\mu}(\theta)(\sigma^2 {\rm
  Id}_T+\tilde{K}^{\mu}(\theta))^{-1},
\end{equation}
and, using the partial sums, noted $K^{\mu,k}_{[n]}$, $k\in V_n$, in (\ref{eq:Kimuinfinite}), to define another sequence $A^{\mu,k}_{[n]}$ which in the limit $n \to \infty$ converge to the coefficients of the Fourier series of $\tilde{A}^\mu$.
We next define a functional $\Gamma_{[n]} = \Gamma_{[n],1} + \Gamma_{[n],2}$, which we use to characterise the Radon-Nikodym derivative of $\Pi^n$ with respect to $R^n$. Let $\mu\in\mathcal{M}_\mS(\mathcal{T}^{\Z})$ and
\begin{equation}\label{eq:Gamma1}
\Gamma_{[n],1}(\mu)=-\frac{1}{2|V_n|}\log\left({\rm det}\left({\rm
  Id}_{|V_n|T}+\frac{1}{\sigma^2} K^{\mu}_{[n]}\right)\right),
\end{equation}
where $K^{\mu}_{[n]}$ is the $(|V_n|T \times |V_n|T)$ covariance matrix of the Gaussian law $\mQ^{\mu}_{[n]}$ defined by the sequence $(K^{\mu,k}_{[n]})_{k\in V_n}$. 

Because of  previous remarks the above expression has a
sense. Taking the limit when $n \to \infty$ does not pose any problem and 
we can define $\Gamma_1(\mu) =\lim_{n\rightarrow\infty}\Gamma_{[n],1}(\mu)$. The following lemma whose proof is straightforward indicates that this is well-defined.
\begin{lemma}\label{lemma:Gamma1}
When $n$ goes to infinity the limit of \eqref{eq:Gamma1} is given by
\begin{equation}\label{eq:Gamma1infinity}
\Gamma_1(\mu)=-\frac{1}{4\pi} \int_{-\pi}^{\pi} \log \left( {\rm det}\left({\rm
  Id}_{T}+\frac{1}{\sigma^2} \tilde{K}^\mu(\theta)\right)\right)\,d\theta
\end{equation}{theo:LDP}
for all $\mu \in \mathcal{M}_{1,S}^+(\T^\Z)$.
\end{lemma}
It also follows easily from previous remarks that
\begin{proposition}\label{prop:Gamma1b}
$\Gamma_{[n],1}$ and $\Gamma_1$ are
bounded below and continuous on  $\mathcal{M}_\mS(\T^\Z)$.
\end{proposition}
The definition of $\Gamma_{[n],2}(\mu)$ is slightly more technical but follows naturally from propositions \ref{prop:RNderiv1} and \ref{prop:radon-nikodym}.
For $\mu \in \mmeas(\T^{\Z})$ let
\begin{equation}\label{eq:Gamma2muN}
\Gamma_{[n],2}(\mu)=\int_{\T_{1,T}^{V_n}}\phi^n(\mu,v)\undt{\mu}^{V_n}_{1,T}(dv)
\end{equation} 
where $\phi^n : \mmeas(\T^{\Z}) \times \T_{1,T}^{V_n} \to \R$ is defined by
\vspace{-0.5cm}
\begin{equation}\label{eq:phiNdeftn}
\phi^n(\mu,v) = \frac{1}{2\sigma^2}\Bigg(\frac{1}{|V_n|}\sum_{j,k \in V_n}
{}^\dagger  (v^j- c^{\mu})A^{\mu,\,k}_{[n]}
   (v^{k+j}-
c^{\mu} )+\\ \frac{2}{|V_n|}\sum_{j \in V_n} \langle c^{\mu}, v^j\rangle -
  \|c^{\mu} \|^2\Bigg).
\end{equation}
$\Gamma_{[n],2}(\mu)$ is finite in the subset $\mathcal{E}_2$ of
$\mM_\mS(\T^{\Z})$ defined in definition \ref{def:E2}. If
$\mu\notin\mathcal{E}_2$, then we set $\Gamma_{[n],2}(\mu) =
\infty$. 

We define
$\Gamma_2(\mu) = \lim_{N\rightarrow\infty}\Gamma_{[n],2}(\mu)$. The following proposition indicates that $\Gamma_2(\mu)$ is well-defined.

\begin{proposition}\label{prop:E2Gamma2}
If the measure $\mu$ is in $\mathcal{E}_2$, i.e. if $\Exp^{\undt{\mu}_{1,T}}[\| v^0 \|^2] < \infty$, then $\Gamma_2(\mu)$ is finite and writes
\begin{multline}\label{eq:Gamma2}
\Gamma_2(\mu)=\frac{1}{2\sigma^2}
\Big(\frac{1}{2\pi}\int_{-\pi}^{\pi} \tilde{A}^\mu(-\theta) :
  \tilde{v}^\mu(d\theta)+
  {}^\dagger  c^\mu
  (\tilde{A}^\mu(0)-{\rm Id}_{T}) c^\mu+\\
  2\Exp^{\undt{\mu}_{1,T}}\left[{}^{t}
   {v^0}({\rm Id}_{T}-\tilde{A}^\mu(0))c^\mu\right]\Big).
\end{multline}
The ``:'' symbol indicates the double
contraction on the indexes. 
\end{proposition}

It is shown in \cite{faugeras-maclaurin:14d} that $\phi^n(\mu,v)$ defined by \eqref{eq:phiNdeftn} is a continuous function of $\mu$ which satisfies
\[
\phi^n(\mu,v) \geq -\beta_2, \quad \beta_2=\frac{T\bar{J}^2}{2\sigma^2 \Lambda^{sum}}(\sigma^2+\theta^2+\Lambda^{asum})
\]
By a standard argument we obtain the following proposition.
\begin{proposition}\label{prop:Gamma2lsc}
$\Gamma_{[n],2}(\mu)$ is lower-semicontinuous.
\end{proposition}
We define $\Gamma_{[n]}(\mu)=\Gamma_{[n],1}(\mu)+\Gamma_{[n],2}(\mu)$. We may conclude from propositions \ref{prop:Gamma1b} and \ref{prop:Gamma2lsc} that $\Gamma_{[n]}$ is lower-semicontinuous hence measurable. 

From these definitions it is relatively easy, and proved in \cite{faugeras-maclaurin:14d},  to show that the measure $Q^{V_n}$ is absolutely continuous with respect to $P^{\otimes V_n}$ with a Radon-Nikodym derivative which can be expressed as a function of the functional $\Gamma_{[n]}$.
\begin{proposition}\label{prop:RNderiv1}
The Radon-Nikodym  derivative
of $Q^{V_n}$ with respect to $P^{\otimes V_n}$ is given by
the following expression.
\begin{equation}\label{eq:dQNdPN}
\frac{dQ^{V_n}}{dP^{\otimes V_n}}(u)=
\Expt{\exp\left(\frac{1}{\sigma^2} \left(\sum_{j \in V_n} \langle \Psi_{1,T}(u^j), G^j \rangle-
\frac{1}{2}\| G^j \|^2\right)\right)},
\end{equation}
for all $u \in V_n$, and the expectation being taken against the $2n+1$ $T$-dimensional Gaussian processes $(G^i)$, $i \in V_n$ given by
\begin{equation*}
G^i_t  =  \sum_{j \in V_n} J_{ij}^{\olivier{N}} f(u^j_{t-1}),
\quad t=1,\cdots,T ,
\end{equation*}
and the function $\Psi$ being defined by \eqref{eq:Psi}.
\end{proposition}
Using standard Gaussian calculus we obtain the following proposition.

\begin{proposition}\label{prop:radon-nikodym}
The Radon-Nikodym derivatives write as 
\begin{align}
\frac{dQ^{V_n}}{dP^{\otimes V_n}}(u) &=\exp(|V_n|\Gamma_{[n]}(\hat{\mu}_n(u)),\label{eq:RN1}\\
\frac{d\Pi^n}{dR^n}(\mu) &= \exp(|V_n|\Gamma_{[n]}(\mu))\label{eq:RN2}.
\end{align}
Here $\mu\in \M_\mS(\T^\Z)$, 
$\Gamma_{[n]}(\mu)=\Gamma_{[n],1}(\mu)+\Gamma_{[n],2}(\mu)$
and the expressions for $\Gamma_{[n],1}$ and
$\Gamma_{[n],2}$ have been defined in equations \eqref{eq:Gamma1} and \eqref{eq:Gamma2muN}. 
\end{proposition}
\vspace{-0.5cm}
\section{The large deviation principle}
\vspace{-0.5cm}
We define the function $H: \mathcal{M}_\mS(\T^\Z) \to [0,+\infty)$ as follows.
\begin{definition}\label{defn:H}
Let $H$ be the function $\mathcal{M}_\mS^+(\T^\Z) \to \R \cup \{+\infty\}$ defined by
\[
H(\mu)=\left\{
\begin{array}{l}
+\infty \quad \text{if} \quad I^{(3)}(\mu,P^\Z)=\infty\\
I^{(3)}(\mu,P^\Z)-\Gamma(\mu) \quad \text{otherwise},
\end{array}
\right.
\]
where $\Gamma=\Gamma_1+\Gamma_2$.
\end{definition}
We finally state the following theorem.
\begin{theorem}\label{theo:LDP}
$\Pi^n$ is governed by a large deviation principle with a good rate
function $H$. 
\end{theorem}

The proof is too long to be reproduced here, see \cite{faugeras-maclaurin:14d}. We only give the general strategy. First we prove the lower bound on the open sets. For the upper bound on the closed sets, we simply avoid it by a) proving that ($\Pi^n$) is exponentially tight which allows us to b) restrict the proof of the upper bound to compact sets. The proof of b) is long and technical. It is partially built upon ideas found in \cite{guionnet:95}.

Note that we have found an analytical form for $H$ through equations \eqref{eq:Gamma1infinity} and \eqref{eq:Gamma2}
\vspace{-0.5cm}





\begin{thebibliography}{10}

\bibitem{cessac:95}
{\sc B.~Cessac}, {\em Increase in complexity in random neural networks},
  Journal de Physique I (France), 5 (1995), pp.~409--432.

\bibitem{cessac-samuelides:07}
{\sc B.~Cessac and M.~Samuelides}, {\em From neuron to neural networks
  dynamics.}, EPJ Special topics: Topics in Dynamical Neural Networks, 142
  (2007), pp.~7--88.

\bibitem{deuschel-stroock-etal:91}
{\sc J.~Deuschel, D.~Stroock, and H.~Zessin}, {\em Microcanonical distributions
  for lattice gases}, Communications in Mathematical Physics, 139 (1991).

\bibitem{faugeras-maclaurin:14d}
{\sc O.~Faugeras and J.~Maclaurin}, {\em {Asymptotic description of neural
  networks with correlated synaptic weights}}, Rapport de recherche RR-8495,
  INRIA, Mar. 2014.

\bibitem{faugeras-maclaurin:14b}
\leavevmode\vrule height 2pt depth -1.6pt width 23pt, {\em Asymptotic
  description of stochastic neural networks. ii - characterization of the limit
  law}, C. R. Acad. Sci. Paris, Ser. I,  (2014).

\bibitem{guionnet:95}
{\sc A.~Guionnet}, {\em Dynamique de Langevin d'un verre de spins}, PhD thesis,
  Universit\'e de Paris Sud, 1995.

\bibitem{moynot:99}
{\sc O.~Moynot}, {\em Etude math\'ematique de la dynamique des r\'eseaux
  neuronaux al\'eatoires r\'ecurrents}, PhD thesis, Universit\'e Paul Sabatier,
  Toulouse, 1999.

\bibitem{moynot-samuelides:02}
{\sc O.~Moynot and M.~Samuelides}, {\em {Large deviations and mean-field theory
  for asymmetric random recurrent neural networks}}, Probability Theory and
  Related Fields, 123 (2002), pp.~41--75.

\bibitem{rassoul-agha:03}
{\sc F.~Rassoul-Agha}, {\em The point of view of the particle on the law of large numbers for random walks in a mixing environment},
The Annals of Applied Probability (2003)

\bibitem{samuelides-cessac:07}
{\sc M.~Samuelides and B.~Cessac}, {\em Random recurrent neural networks},
  European Physical Journal - Special Topics, 142 (2007), pp.~7--88.

\bibitem{sompolinsky-crisanti-etal:88}
{\sc H.~Sompolinsky, A.~Crisanti, and H.~Sommers}, {\em {Chaos in Random Neural
  Networks}}, Physical Review Letters, 61 (1988), pp.~259--262.

\bibitem{sznitman-zerner:99}{\sc A.~ Sznitman and M.~Zerner} {\em {A Law of Large Numbers for Random Walks in Random Environment}}, The Annals of Probability, 27 (1999), pp.~1851-1869.


\end{thebibliography}

\end{document}